\documentclass[11pt]{article}

\usepackage{latexsym}
\usepackage{amssymb}

\def\C{{\mathbb C}}
\def\N{{\mathbb N}}
\def\Z{{\mathbb Z}}

\def\Mat{{\rm Mat}}

\def\Ad{{\rm Ad}\,}

\def\Tr{{\rm Tr}}
\def\supp{{\rm supp}\,}

\newtheorem{theorem}{Theorem}
\newtheorem{lemma}[theorem]{Lemma}
\newtheorem{corollary}[theorem]{Corollary}
\newtheorem{proposition}[theorem]{Proposition}

\newenvironment{verif}{\noindent{\it Proof.}}
{\newline\mbox{\ }\hfill\rule{0.5em}{0.5em}\smallskip}

\begin{document}

\title{\bf Entropy in type I algebras}

\author{
Sergey Neshveyev$^{1)}$ \\
{\it\small B. Verkin Institute for Low
Temperature Physics and Engineering,}\\
{\it\small 47, Lenin Ave., 310164, Kharkov, Ukraine}\\
\\
and\\
\\
Erling St{\o}rmer$^{2)}$ \\
{\it\small Department of Mathematics, University of Oslo,}\\
{\it\small P.O. Box 1053, Blindern, 0316 Oslo, Norway}}

\date{}

\footnotetext[1]{ Partially supported by NATO grant SA
(PST.CLG.976206)5273.}
\footnotetext[2]{ Partially supported by the Norwegian Research Council.}

\maketitle

\begin{abstract}
It is shown that if $(M,\phi,\alpha)$ is a
W$^*$-dynamical system with $M$ a type~I von~Neumann algebra then
the entropy of $\alpha$ w.r.t. $\phi$ equals the entropy of the
restriction of $\alpha$ to the center of $M$. If furthermore
$(N,\psi,\beta)$ is a W$^*$-dynamical system with $N$ injective
then
$h_{\phi\otimes\psi}(\alpha\otimes\beta)
=h_{\phi}(\alpha)+h_{\psi}(\beta)$.
\end{abstract}


\section{Introduction}

In the theory of non-commutative entropy the attention has almost
exclusively been concentrated on non type~I algebras. We shall in
the present paper remedy this situation by proving the basic
facts on entropy of automorphisms of type~I C$^*$- and
von Neumann-algebras. The results are as nice as one can hope. The
CNT-entropy of an automorphism of a von~Neumann algebra of type~I
with respect to an invariant normal state is the classical entropy of the
restriction of the automorphism to the center of the algebra. If
one factor of a tensor product of two von~Neumann algebras is of
type~I and the other injective, then the entropy of a tensor product
automorphism with respect to an invariant product state is the sum of the
entropies. The results have obvious corollaries to type~I C$^*$-algebras.
The  main idea behind our proofs is the use of conditional expectations
of finite index, as employed in~\cite{GN}.

We shall use the notation $h_\phi(\alpha)$ for the CNT-entropy of
a C$^*$-dynamical system as defined by Connes, Narnhofer and
Thirring in~\cite{CNT}, and $h'_{\phi}(\alpha)$ for the ST-entropy
defined by Sauvageot and Thouvenot in~\cite{ST}.

\bigskip

\section{Main results}

We first prove a general result for the Sauvageot-Thouvenot
entropy for the restriction of an automorphism to a globally
invariant C$^*$-subalgebra of finite index. Recall the definition
of ST-entropy and its connection with CNT-entropy.

A stationary coupling of a C$^*$-dynamical system
$(A,\phi,\alpha)$ with a commutative system $(C,\mu,\beta)$ is an
$\alpha\otimes\beta$-invariant state $\lambda$ on $A\otimes C$
such that $\lambda|_A=\phi$ and $\lambda|_C=\mu$. Given such a
coupling and a finite-dimensional subalgebra $P$ of $C$ with atoms
$p_1,\ldots,p_n$, consider the quantity
$$
H_\mu(P|P^-)-H_\mu(P)+\sum^n_{i=1}\mu(p_i)S(\phi,\phi_i),
$$
where $\phi_i(a)={1\over\mu(p_i)}\lambda(a\otimes p_i)$. By
definition, the ST-entropy $h'_\phi(\alpha)$ of the system
$(A,\phi,\alpha)$ is the supremum of these quantities.

By \cite[Proposition 4.1]{ST}, ST-entropy coincides with
CNT-entropy for nuclear C$^*$-algebras. In fact, the proof of the
inequality $h_\phi(\alpha)\le h'_\phi(\alpha)$ does not use any
assumptions on the algebra. On the other hand, given a coupling
$\lambda$ and an algebra $P$ as above, for each $m\in\N$ we can
form the decomposition
$$
\phi=\sum^n_{i_1,\ldots,i_m=1}\phi_{i_1\ldots i_m},\ \
\phi_{i_1\ldots i_m}(a)=\lambda(a\otimes
p_{i_1}\beta(p_{i_2})\ldots\beta^{m-1}(p_{i_m})).
$$
If $\gamma$ is a unital completely positive mapping of a
finite-dimensional C$^*$-algebra into $A$, we can use these
decompositions in computing the mutual entropy
$H_\phi(\gamma,\alpha\circ\gamma,\ldots,\alpha^{m-1}\circ\gamma)$
\cite{CNT}. Indeed, since the atoms in $\beta^j(P)$ are
$\beta^j(p_1),\ldots,\beta^j(p_n)$ we have by \cite[III.3]{CNT}

$\displaystyle
\lefteqn{H_\phi(\gamma,\alpha\circ\gamma,\ldots,\alpha^{m-1}\circ\gamma)
\geq S\Bigl(\mu\Big|\bigvee_0^{m-1}\beta^j(P)\Bigr)-\sum_{j=0}^{m-1}
S\Bigl(\mu\Big| \beta^j(P)\Bigr)}$
$$\hspace{4cm}
+\sum_j\sum_i \mu(\beta^j(p_i))
S\Bigl(\phi\circ \alpha^j\circ \gamma,\frac
{\lambda((\alpha^j\circ\gamma)(\cdot)\otimes\beta^j(p_i))}
{\mu(\beta^j(p_i))}\Bigr).
$$
Hence by invariance of $\phi$, $\mu$ and $\lambda$ with respect to $\alpha$,
$\beta$ and $\alpha\otimes\beta$ respectively
$$
\frac{1}{m}H_\phi(\gamma,\alpha\circ\gamma,\ldots,
\alpha^{m-1}\circ\gamma)\geq \frac{1}{m}H_\mu\Bigl(\bigvee_0^{m-1}
\beta^j(P)\Bigr) -H_\mu(P)+\sum_i
\mu(p_i)S(\phi\circ\gamma,\phi_i\circ\gamma).
$$
It follows that
$$
h_\phi(\alpha)\ge
H_\mu(P|P^-)-H_\mu(P)+\sum^n_{i=1}\mu(p_i)
S(\phi\circ\gamma,\phi_i\circ\gamma).
$$
Thus what is really necessary for the coincidence of the
entropies, is the existence of a net of unital completely positive
mappings $\gamma_i$ of finite-dimensional C$^*$-algebras into $A$
such that
$S(\phi,\psi)=\lim_iS(\phi\circ\gamma_i,\psi\circ\gamma_i)$ for
any positive linear functional $\psi$ on $A$, $\psi\le\phi$. In
particular, $h_\phi(\alpha)=h'_\phi(\alpha)$ if $A$ is an
injective von Neumann algebra and $\phi$ is a normal state on it.

\begin{proposition} \label{1}
Let $(A,\phi,\alpha)$ be a unital
C$^*$-dynamical system. Let $B\subset A$ be an $\alpha$-invariant
C$^*$-subalgebra (with $1\in B$). Suppose there exists a
conditional expectation $E\colon A\to B$ such that
$E\circ\alpha=\alpha\circ E$, $\phi\circ E=\phi$ and $E(x)\ge cx$
for all $x\in A^+$ for some $c>0$. Then
$h'_{\phi}(\alpha)=h'_{\phi}(\alpha|_B)$.
\end{proposition}

\begin{verif}
Let $(C,\mu,\beta)$ be a C$^*$-dynamical system with $C$ abelian. Using
$E$
we can lift any stationary coupling on $B\otimes C$ to a
stationary coupling on $A\otimes C$. This, together with the property of
monotonicity of relative entropy, shows that $h'_{\phi}(\alpha)\ge
h'_{\phi}(\alpha|_B)$.

Conversely, suppose $\lambda$ is a stationary coupling of
$(A,\phi,\alpha)$ with $(C,\mu,\beta)$, $P$ a finite-dimen\-sional
subalgebra of $C$  with atoms $p_1,\ldots,p_n$, and
$\phi_i(a)={1\over\mu(p_i)}\lambda(a\otimes p_i)$ for $a\in A$.
Since $\phi_i\le{1\over\mu(p_i)}\phi$, $\phi_i$ is normal in the
GNS-representation of $\phi$. Since $E$ is $\phi$-invariant, it
extends to a normal conditional expectation of the closure of $A$
in the GNS-representation onto the closure of $B$. Thus we can
apply~\cite[Theorem 5.15]{OP} to $\phi$ and $\phi_i$, and (as in
the proof of Lemma~1.5 in~\cite{GN}) get
$$
\sum^n_{i=1}\mu(p_i)S(\phi,\phi_i)
=\sum^n_{i=1}\mu(p_i)(S(\phi|_B,\phi_i|_B)+S(\phi_i\circ E,\phi_i))
\le\sum^n_{i=1}\mu(p_i)S(\phi|_B,\phi_i|_B)-\log c.
$$
It follows that $h'_{\phi}(\alpha)\le h'_{\phi}(\alpha|_B)-\log c$. Then
for each $m\in\N$
$$
h'_\phi(\alpha)={1\over m}h'_\phi(\alpha^m)
\le{1\over m}h'_\phi(\alpha^m|_B)-{1\over m}\log c
=h'_\phi(\alpha|_B)-{1\over m}\log c.
$$
Thus $h'_\phi(\alpha)\le h'_\phi(\alpha|_B)$.
\end{verif}

\begin{corollary} \label{2}
If in the above proposition $A$ and $B$ are injective von~Neumann algebras
and $\phi$ is normal then $h_\phi(\alpha)=h_\phi(\alpha|_B)$.
\end{corollary}

To prove our main result we need also two simple lemmas. The first lemma
is more or less well-known.

\begin{lemma} \label{3}
Let $(M,\phi,\alpha)$ be a W$^*$-dynamical system. Then

(i) if $p$ is an
$\alpha$-invariant projection in $M$ such that $\supp\phi\le p$, then
$h_\phi(\alpha)=h_\phi(\alpha|_{M_p})$;

(ii) if $\{p_i\}_{i\in I}$ is a set of mutually orthogonal
$\alpha$-invariant central projections in $M$, $\sum_ip_i=1$, then
$$
h_\phi(\alpha)=\sum_i\phi(p_i)h_{\phi_i}(\alpha_i),
$$
where $\phi_i={1\over\phi(p_i)}\phi$ is the normalized restriction of
$\phi$ to $Mp_i$, and
$\alpha_i=\alpha|_{Mp_i}$.
\end{lemma}

\begin{verif}
(i) easily follows from the definitions; (ii) follows
from~\cite[VII.5(iii)]{CNT}, (i) and~\cite[Lemma 3.3]{SV} applied
to the subalgebras
$M(p_{i_1}+\ldots+p_{i_n})+\C(1-p_{i_1}-\ldots-p_{i_n})$.
\end{verif}

The proof of the following lemma is left to the reader.

\begin{lemma} \label{3a}
Let $T$ be an automorphism of a probability space $(X,\mu)$, $f\in
L^\infty(X,\mu)$ a $T$-invariant function such that $f\ge0$
and $\int_Xf\,d\mu=1$. Let $\mu_f$ be the measure on $X$ such that
$d\mu_f/d\mu=f$. Then $h_{\mu_f}(T)\le||f||_\infty h_\mu(T)$.
\end{lemma}

\begin{theorem} \label{4}
Let $(M,\phi,\alpha)$ be a
W$^*$-dynamical system with $M$ a von~Neumann algebra of type~I.
Let $Z$ denote the center of $M$. Then
$h_\phi(\alpha)=h_\phi(\alpha|_Z)$.
\end{theorem}

\begin{verif}
By Lemma~\ref{3}(i) we may suppose that $\phi$ is
faithful. Then $M$ is a direct sum of homogeneous algebras of type
I$_n$, $n\in\N\cup\{\infty\}$. By Lemma~\ref{3}(ii) we may assume
that $M$ is homogeneous of type I$_n$. We first assume that
$n\in\N$. Then $Z=L^\infty(X,\mu)$, where $(X,\mu)$ is a
probability space and $\phi|_Z=\mu$. Thus
$$
M\cong
Z\otimes\Mat_n(\C)=L^\infty(X,\Mat_n(\C)),\ \
\phi=\int^\oplus_X\phi_xd\mu(x),
$$
where $\phi_x=\Tr(\cdot\,Q_x)$
is a state on $\Mat_n(\C)$, $\Tr$ the canonical trace on
$\Mat_n(\C)$. We first assume $Q_x\ge c>0$ for all $x$.

If $s\in M^+$, $s$ is a function in $L^\infty(X,\Mat_n(\C))$. Define the
$\phi$-preserving conditional expectation $E\colon M\to Z$ by
$E(s)(x)=\phi_x(s(x))$. Then
$$
E(s)(x)=\Tr(s(x)Q_x)\ge c\Tr(s(x))\ge cs(x),
$$
so $E(s)\ge cs$, and it follows from Corollary~\ref{2} that
$h_\phi(\alpha)=h_\phi(\alpha|_Z)$.

If there is no $c>0$ such that $Q_x\ge c$ for all $x$, let
$X_c=\{x\in X\,|\, Q_x\ge c\}$, ($c>0$),
$$
N_c=L^\infty(X_c,\Mat_n(\C))\ \ \hbox{and}\ \
M_c=N_c+\C\chi_{X\backslash X_c},
$$
where $\chi_{X\backslash X_c}$ is the characteristic function of
$X\backslash X_c$. Since $\phi$ is $\alpha$-invariant so is $M_c$, so by
the above argument and Lemma~\ref{3}, letting
$\phi_c={1\over\mu(X_c)}\phi|_{N_c}$ and
$\mu_c={1\over\mu(X_c)}\mu|_{X_c}$, we obtain
$$
h_\phi(\alpha|_{M_c})=\mu(X_c)h_{\phi_c}(\alpha|_{N_c})
=\mu(X_c)h_{\mu_c}(T|_{X_c})
\le h_\mu(T),
$$
where $T$ is the automorphism of $(X,\mu)$ induced by $\alpha$.
Letting $c\to0$ and using~\cite[Lemma~3.3]{SV} we obtain the Theorem when
$M$  is finite.

If $M$ is homogeneous of type I$_\infty$, we have
$M\cong L^\infty(X,\mu)\otimes B(H)$, where $H$ is a separable Hilbert
space. Let $\Tr$ denotes the canonical trace on $B(H)$. Write again
$$
\phi=\int^\oplus_X\phi_xd\mu(x),\ \
\phi_x=\Tr(\cdot\,Q_x),
$$
and let $E_x(U)$ denote the spectral projection of $Q_x$ corresponding
to a Borel set $U$. Let $P_c\in M=L^\infty(X,B(H))$ be the projection
defined by $P_c(x)=E_x([c,+\infty))$, where $c>0$. Then $P_c$ is an
$\alpha$-invariant finite projection. Let
$$
M_c=P_cMP_c+\C(1-P_c).
$$
Then $M_c$ is a finite type I von~Neumann algebra. Its center is
isomorphic
to $L^\infty(X_c,\mu_c)\oplus\C$, and the restriction of $\phi$ to it is
$\phi(P_c)\mu_c\oplus\phi(1-P_c)$, where $X_c=\{x\in X\,|\,P_c(x)\ne0\}$
and
$$
\int_{X_c}f(x)d\mu_c(x)
={1\over\phi(P_c)}\int_{X_c}f(x)\phi_x(P_c(x))d\mu(x).
$$
So we can apply the first part of the proof to $M_c$. Since
$d\mu_c/d\mu\le{1\over\phi(P_c)}$, applying Lemma~\ref{3a} we get
$$
h_\phi(\alpha|_{M_c})=\phi(P_c)h_{\mu_c}(T|_{X_c})\le h_\mu(T).
$$
Now letting $c\to0$ we conclude that $h_\phi(\alpha)=h_\mu(T)$.
\end{verif}

It should be remarked that in a special case the above theorem was
proved in~\cite[Proposition~2.4]{GS1}.
\medskip

If $A$ is a C$^*$-algebra and $\phi$ a state on $A$, the central measure
$\mu_\phi$ of $\phi$ is the measure on the spectrum $\hat A$ of $A$
defined
by $\mu_\phi(F)=\phi(\chi_F)$, where $\phi$ is regarded as a normal state
on $A''$, see~\cite[4.7.5]{P}. Thus by Theorem~\ref{4} and~\cite[4.7.6]{P}
we have the following

\begin{corollary} \label{5}
Let $(A,\phi,\alpha)$ be a C$^*$-dynamical system with $A$ a separable
unital type I C$^*$-algebra. Then
$h_\phi(\alpha)=h_{\mu_\phi}(\hat\alpha)$,
where $\hat\alpha$ is the automorphism of the measure space
$(\hat A,\mu_\phi)$ induced by $\alpha$.
\end{corollary}

Since inner automorphisms act trivially on the center we have

\begin{corollary} \label{6}
If $(M,\phi,\alpha)$ is a W$^*$-dynamical system with $M$ of type I and
$\alpha$ an inner automorphism then $h_\phi(\alpha)=0$.
\end{corollary}

Note that in the finite case the above corollary also follows from a
result of N.~Brown~\cite[Lemma 2.2]{Br}.
\medskip

The next result was shown in~\cite{S} when $\phi$ is a trace.

\begin{corollary} \label{7}
Let $R$ denote the hyperfinite II$_1$-factor. Let $A$ be a Cartan
subalgebra of $R$ and $u$ a unitary operator in $A$. If $\phi$ is a normal
state such that $u$ belongs to the centralizer of $\phi$ then
$h_\phi(\Ad u)=0$.
\end{corollary}

\begin{verif}
As in~\cite{S}, it follows from~\cite{CFW} that there exists an increasing
sequence of full matrix algebras $N_1\subset N_2\subset\ldots$ with union
weakly dense in $R$ such that $A\cong A_n\otimes B_n$, where
$A_n=N_n\cap A$ and $B_n=(N_n'\cap R)\cap A$ for all $n\in\N$. Let
$M_n=N_n\otimes B_n$.
Then $M_n$ is of type I and contains~$u$. Hence
$h_\phi(\Ad u|_{M_n})=0$. Since $(\cup_nM_n)^-=R$, $h_\phi(\Ad u)=0$
by~\cite[Lemma~3.3]{SV}.
\end{verif}

If $(A,\phi,\alpha)$ and $(B,\psi,\beta)$ are C$^*$-dynamical
systems we always have
$$
h_{\phi\otimes\psi}(\alpha\otimes\beta)
\ge h_\phi(\alpha)+h_\psi(\beta),
$$
see \cite[Lemma~3.4]{SV}. The equality
does not always hold, see~\cite{NST} or~\cite{Sa}. However, we
have

\begin{theorem} \label{8}
Let $(A,\phi,\alpha)$ and
$(B,\psi,\beta)$ be W$^*$-dynamical systems. Suppose that $A$ is
of type~I, and $B$ is injective. Then
$$
h_{\phi\otimes\psi}(\alpha\otimes\beta)=
h_\phi(\alpha)+h_\psi(\beta).
$$
\end{theorem}

\begin{verif}
We shall rather prove that
$h_{\phi\otimes\psi}(\alpha\otimes\beta)
 =h_\phi(\alpha|_{Z(A)})+h_\psi(\beta)$. For this
it suffices to consider the case when $A$ is abelian; the general case
will follow by the same arguments as in the proof of Theorem~\ref{4}.
(Note that the mapping $x\mapsto\Tr(x)-x$ on
$\Mat_n(\C)$ is not completely positive, but the mapping
$x\mapsto\Tr(x)-{1\over n}x$ is by the Pimsner-Popa
inequality. Thus replacing $M$ with $M\otimes B$ and $Z$ with
$Z\otimes B$ in the proof of
Theorem~\ref{4} we have to replace the inequality
$E(s)\ge cs$ in the proof with $E(s)\ge{c\over n}s$.)

So suppose that $A$ is abelian. It is clear that it
suffices to prove that if $A_1,\ldots,A_n$ are finite-dimensional
subalgebras of $A$, and $B_1,\ldots,B_n$ are finite-dimensional
subalgebras of $B$, then
$$
H_{\phi\otimes\psi}(A_1\otimes B_1,\ldots,A_n\otimes
B_n)=H_\phi(A_1,\ldots,A_n)+H_\psi(B_1,\ldots,B_n).
$$
We always
have the inequality "$\ge$", \cite[Lemma~3.4]{SV}. To prove the opposite
inequality consider a decomposition
$$
\phi\otimes\psi=\sum_{i_1,\ldots,i_n}\omega_{i_1\ldots i_n}.
$$
Let $H_{\{\phi\otimes\psi=\sum\omega_{i_1\ldots i_n}\}}(A_1\otimes
B_1,\ldots,A_n\otimes B_n)$ be the entropy of the corresponding
abelian model, so \newline $\displaystyle
H_{\{\phi\otimes\psi=\sum\omega_{i_1\ldots i_n}\}}(A_1\otimes
B_1,\ldots,A_n\otimes B_n)=$
$$
=\sum_{i_1,\ldots,i_n}\eta\omega_{i_1\ldots
i_n}(1)+\sum^n_{k=1}\sum_iS\left(\phi\otimes\psi|_{A_k\otimes
B_k},\sum_{i_k=i}\omega_{i_1\ldots i_n}|_{A_k\otimes B_k}\right).
$$
Set $C=\vee^n_{k=1}A_k$. Let $p_1,\ldots,p_r$ be those atoms
$p$ of $C$ for which $\phi(p)>0$. Define positive linear
functionals $\psi_{m,i_1\ldots i_n}$ on $B$,
$$
\psi_{m,i_1\ldots i_n}(b)={\omega_{i_1\ldots i_n}(p_m\otimes b)
\over\phi(p_m)}.
$$
Let also $\phi_m$ be the linear functional
on $C$ defined by the equality $\phi_m(a)=\phi(ap_m)$. Then
$$
\omega_{i_1\ldots i_n}=\sum^r_{m=1}\phi_m\otimes
\psi_{m,i_1\ldots i_n}\ \ \hbox{on}\ \ C\otimes B,
$$
and
$$
\psi=\sum_{i_1,\ldots,i_n}\psi_{m,i_1\ldots i_n}\ \ \hbox{for}\
\ m=1,\ldots,r.
$$
Since the supports of the states $\phi_m$ are mutually orthogonal minimal
projections in $C$, we have
$\displaystyle
\sum^n_{k=1}\sum_iS\left(\phi\otimes\psi|_{A_k\otimes
B_k},\sum_{i_k=i}\omega_{i_1\ldots i_n}|_{A_k\otimes
B_k}\right)\le $
\begin{eqnarray*}
&\le&\sum^n_{k=1}\sum_iS\left(\phi\otimes\psi|_{C\otimes B_k},
        \sum_{i_k=i}\omega_{i_1\ldots i_n}|_{C\otimes B_k}\right)\\
&=&\sum^n_{k=1}\sum_iS\left(\phi\otimes\psi|_{C\otimes B_k},
    \sum^r_{m=1}\phi_m\otimes\left(\sum_{i_k=i}\psi_{m,i_1\ldots i_n}
      \right)|_{C\otimes B_k}\right)\\
&=&\sum^n_{k=1}\sum_i\sum^r_{m=1}\phi(p_m)S\left(\psi|_{B_k},
      \sum_{i_k=i}\psi_{m,i_1\ldots i_n}|_{B_k}\right).
\end{eqnarray*}
If $a_i\geq 0$ and $\sum\limits_i a_i\leq 1$ then $\eta(\sum\limits_i
a_i)\leq \sum\limits_i \eta(a_i)$. Hence we have
\begin{eqnarray*}
\sum_{i_1,\ldots,i_n}\eta\omega_{i_1\ldots i_n}(1)
&\le&\sum^r_{m=1}\sum_{i_1,\ldots,i_n}\eta
      (\phi_m\otimes\psi_{m,i_1\ldots i_n})(1)\\
&=&\sum^r_{m=1}\eta\phi(p_m)\sum_{i_1,\ldots,i_n}
     \psi_{m,i_1\ldots i_n}(1)
     +\sum^r_{m=1}\phi(p_m)\sum_{i_1,\ldots,i_n}\eta\psi_{m,i_1\ldots
    i_n}(1)\\
&=&\sum^r_{m=1}\eta\phi(p_m)+\sum^r_{m=1}
     \phi(p_m)\sum_{i_1,\ldots,i_n}\eta\psi_{m,i_1\ldots i_n}(1).
\end{eqnarray*}
Thus
\begin{eqnarray*}
\lefteqn{H_{\{\phi\otimes\psi=\sum\omega_{i_1\ldots i_n}\}}
     (A_1\otimes B_1,\ldots,A_n\otimes B_n)\le} \\
&&\qquad \le\sum^r_{m=1}\eta\phi(p_m)+
    \sum^r_{m=1}\phi(p_m)
      H_{\{\psi=\sum\psi_{m,i_1\ldots i_n}\}}(B_1,\ldots,B_n).
\end{eqnarray*}

Since $\sum_m\eta\phi(p_m)=H_\phi(C)=H_\phi(A_1,\ldots,A_n)$, we conclude
that
$$
H_{\phi\otimes\psi}(A_1\otimes B_1,\ldots,A_n\otimes B_n)
\le H_\phi(A_1,\ldots,A_n)+H_\psi(B_1,\ldots,B_n),
$$
completing the proof of the Theorem.
\end{verif}

\end{document}